\documentclass[12pt,twoside]{amsart}
\usepackage{amssymb}
\usepackage{amscd}
\usepackage{mathrsfs}
\usepackage[abbrev,alphabetic]{amsrefs}
\usepackage{hyperref}
\usepackage{tikz-cd}
\usepackage{comment}
\usetikzlibrary{arrows}

\usepackage[margin=1.25in]{geometry}

\title[Kawamata-Viehweg vanishing for toric varieties]
{Kawamata-Viehweg vanishing for toric varieties} 
\author{Hiromu Tanaka} 
\address{Graduate School of Mathematical Sciences, 
The University of Tokyo, 
3-8-1 Komaba, Meguro-ku, Tokyo 153-8914, JAPAN} 
\email{tanaka@ms.u-tokyo.ac.jp}

\newcommand{\Proj}[0]{{\operatorname{Proj}}}
\newcommand{\Spec}[0]{{\operatorname{Spec}}}

\newcommand{\Hom}[0]{{\operatorname{Hom}}}

\newcommand{\Supp}[0]{{\operatorname{Supp}}}

\newcommand{\Ex}[0]{{\operatorname{Ex}}}

\newcommand{\NE}[0]{{\operatorname{NE}}}
\newcommand{\Relint}[0]{{\operatorname{Relint}}}
\newtheorem{thm}{Theorem}[section]
\newtheorem{lem}[thm]{Lemma}
\newtheorem{cor}[thm]{Corollary}
\newtheorem{prop}[thm]{Proposition}

\theoremstyle{definition}

\newtheorem{dfn}[thm]{Definition}

\newtheorem{rem}[thm]{Remark}
\newtheorem{ex}[thm]{Example}

\newtheorem{step}{Step}
\newtheorem{nothing}[thm]{}

\makeatletter
  
  \@addtoreset{equation}{thm}
  \makeatother

\newcommand{\cred}{\color{black}}

\newcommand{\MO}{\mathcal{O}}

\newcommand{\R}{\mathbb{R}}
\newcommand{\Q}{\mathbb{Q}}
\newcommand{\Z}{\mathbb{Z}}

\newcommand{\G}{\mathbb{G}}

\begin{document}

\maketitle

\begin{abstract}
Given a boundary divisor $B$ on a projective toric variety $X$ such that $(X, B)$ is klt, we establish the Kawamata-Viehweg vanishing theorem for $(X, B)$. 
\end{abstract}

\tableofcontents


\section{Introduction}


Many vanishing theorems have been established for toric varieties \cite{Dem70}, \cite{BTLM97}, \cite{Mus02}, \cite{Fuj07} (cf. \cite[Section 9]{CLS11}). 
Amongst them, a typical result is the Danilov vanishing theorem, which states that 
given  a smooth projective toric variety $X$ and an ample Cartier divisor $A$, 
the equation 
\[
H^i(X, \Omega^a_X \otimes \MO_X(A))=0
\]
holds for any $i>0$ and $a \geq 0$ (cf. \cite[Theorem 9.3.1]{CLS11}, \cite[Corollary 1.3]{Fuj07}). 
In particular, smooth proper toric varieties satisfy the Kodaira vanishing even in positive characteristic. 
As a generalisation of the Kodaira vanishing theorem, 
also the Kawamata--Viehweg vanishing theorem holds for a klt pair $(X, B)$, 
where $X$ is a proper toric variety and $B$ is a torus invariant effective $\Q$-divisor \cite[Corollary 1.7]{Fuj07}.  
The main purpose of this note is to drop the assumption that $B$ is torus invariant.

\begin{thm}[Theorem \ref{t-main-full}]\label{t-intro-main1}
Let $k$ be a field. 
Let $\alpha: X \to S$ be a projective $k$-morphism of 
schemes of finite type over $k$, where 
$X$ is a normal toric variety over $k$. 
Let $D$ be a $\Q$-Cartier $\Z$-divisor on $X$. 
Assume that one of the following conditions holds. 
\begin{enumerate}
\item $D \sim_{\R} K_X+B$ for some effective $\R$-Cartier $\R$-divisor $B$ on $X$ 
such that $(X, B)$ is klt and $B$ is $\alpha$-big. 
\item $D - (K_X+B)$ is $\alpha$-nef and $\alpha$-big for some effective $\R$-Cartier $\R$-divisor $B$ on $X$ 
such that $(X, B)$ is klt. 
\end{enumerate}
Then $R^i\alpha_*\MO_X(D)=0$ for any $i>0$. 
\end{thm}


\subsection{Sketch of the proof}

We give a sketch of the proof of Theorem \ref{t-intro-main1}. 
For simplicity, we only treat the case when $S = \Spec\,k$. 
Taking a toric small $\Q$-factorialisation, 
we may assume that $X$ is $\Q$-factorial (cf. the proof of Theorem \ref{t-main-full}). 
Then the assumptions (1) and (2) are known to be equivalent if $k$ is of characteristic zero. 
Although we need some argument, 
we can show that (2) implies (1) under the assumption that $X$ is $\Q$-factorial 
(cf. the proof of Theorem \ref{t-main-full}). 
In what follows, we assume that (1) holds. 
The main idea is to run a $D$-MMP: 
\[
X =: X_0 \overset{f_0}{\dashrightarrow} X_1 \overset{f_1}{\dashrightarrow} \cdots \overset{f_{N-1}}{\dashrightarrow} X_N
\]
\[
D=:D_0, \qquad (f_n)_*D_n =:D_{n+1}.
\]
Such a strategy is used in \cite{CTW17} to establish the Kawamata--Viehweg vanishing theorem for 
klt del Pezzo surfaces in larger characteristic. 
It is enough to prove the following assertions. 
\begin{enumerate}
\renewcommand{\labelenumi}{(\alph{enumi})}
\item $H^i(X_n, \MO_{X_n}(D_n)) \simeq H^i(X_{n+1}, \MO_{X_{n+1}}(D_{n+1}))$ holds 
for any $i>0$ and $n \in \{0, ..., N-1\}$. 
\item $H^i(X, \MO_{X}(D)) =0$ holds for any $i>0$ 
under the assumption that 
either $D$ is nef or $X$ has a $D$-Mori fibre space $f: X \to Y$. 
\end{enumerate}

(a) 
Recall that $f_n:X_n \dashrightarrow X_{n+1}$ is either a $D_n$-divisorial contraction or a $D_n$-flip. 
If $f_n$ is a $D_n$-divisorial contraction, then 
it is not difficult to check that we can apply the relative  Kawamata--Viehweg vanishing theorem for $f_n$ by using the fact that the $f_n$-exceptional divisor is torus invariant. 
We then obtain the required isomorphism $H^i(X_n, \MO_{X_n}(D_n)) \simeq H^i(X_{n+1}, \MO_{X_{n+1}}(D_{n+1}))$. 

Assume that $f_n$ is a flip. 
In this case, it does not seem to be easy to compare 
$H^i(X_n, \MO_{X_n}(D_n))$ and $H^i(X_{n+1}, \MO_{X_{n+1}}(D_{n+1}))$ via $Z$, where 
$\varphi : X_n \to Z$ denotes the corresponding flipping contraction. 
We use the following diagram (\ref{intro-e1}), where $Y \to X_n$ and $Y \to X_{n+1}$ are suitable divisorial contractions. 
Such a diagram always exists as we consider a toric MMP (Proposition \ref{p-flip-blowup}). 
Then the situation is quite similar to the divisorial contraction, although we need further delicate computation. 
For the details on this computation, see Step \ref{s4-toric-Qfac} of the proof of Theorem \ref{t-toric-Qfac}.
\begin{equation}\label{intro-e1}
\begin{tikzcd}
	& Y \\
	{X_n} && {X_{n+1}} \\
	& Z
	\arrow["\psi"', from=1-2, to=2-1]
	\arrow["\varphi"', from=2-1, to=3-2]
	\arrow["{\psi'}", from=1-2, to=2-3]
	\arrow["{\varphi'}", from=2-3, to=3-2]
	\arrow["{f_n}", dashed, from=2-1, to=2-3]
\end{tikzcd}
\end{equation}

(b) 
If $D$ is nef, it is well known that the assertion (b) is true (cf. Proposition \ref{p-GFR-easy}(2)). 
Therefore, the problem is reduced to the case when $X$ has a $D$-Mori fibre space $f: X \to Y$. 
Note that $D-K_X$ is $f$-ample and $f$ is a toric morphism to a normal projective toric variety $Y$. 
Then we have $R^if_*\MO_X(D)=0$ for any $i >0$ (cf. Proposition \ref{p-GFR-easy}(1)). 
Since $-D$ is $f$-big and $\dim X > \dim Y$,  we also have $f_*\MO_X(D)=0$. 
Therefore, it holds that $Rf_*\MO_X(D)=0$. 
For the  morphisms $\alpha : X \xrightarrow{f} Y \xrightarrow{\beta} \Spec\,k$, we have 
the following isomorphism in the derived category: 
\[
R\alpha_* \MO_X(D) \simeq R\beta_* Rf_*\MO_X(D) =0. 
\]
This implies that $H^i(X, \MO_X(D))=0$ for any $i \geq 0$. 

\medskip

Some of the above arguments work even for globally $F$-regular varieties. 
Indeed, the main theorem holds for globally $F$-regular surfaces by a similar and easier proof. 

\begin{thm}[Theorem \ref{t-GFR-surface}]\label{intro-GFR-surface}
Assume that $k$ is of characteristic $p>0$. 
Let $\alpha :X \to S$ be a proper morphism 
from a normal surface $X$ to a scheme $S$ of finite type over $k$ 
such that  $X$ is globally $F$-regular over $S$. 
Let $D$ be a 
$\Z$-divisor on $X$ such that 
$D -(K_X+B)$ is $\alpha$-nef and $\alpha$-big 
for some $\R$-divisor $B$ on $X$ whose coefficients are contained in $[0, 1)$. 
Then $R^i\alpha_*\MO_X(D)=0$ for any $i > 0$. 
\end{thm}

\begin{rem}
\begin{enumerate}
\item 
The author does not know whether Theorem \ref{t-intro-main1} holds after replacing \lq\lq toric" by \lq\lq globally $F$-regular". 
\item 
Theorem \ref{t-intro-main1} is known if $\dim X=2$ \cite[Theorem A]{WX}. 
Both Theorem \ref{t-intro-main1} and Theorem \ref{intro-GFR-surface} can be considered as generalisations of \cite[Theorem A]{WX}. 
\end{enumerate}
\end{rem}



\vspace{5mm}

\textbf{Acknowledgements:} 
The author would like to thank 
Professors  Michel Brion, Osamu Fujino, and Jakub Witaszek 
for useful comments, constructive suggestions, and answering questions. 
{\cred The author also thanks the referee for reading the manuscript carefully and for suggesting several  improvements.} 
The author was funded by JSPS KAKENHI Grant numbers 
JP18K13386, JP22H01112, and JP23K03028. 

\section{Preliminaries}

\subsection{Notation}

In this subsection, we summarise notation and terminologies used in this paper. 

\begin{enumerate}
\item 
We will freely use the notation and terminology in \cite{Har77} and \cite{Kol13}. 
\item 
Throughout this paper, we work over a field $k$. 
\item 
We say that $X$ is a {\em variety} (over $k$) 
if $X$ is an integral scheme that is separated and of finite type over $k$. 
We say that $X$ is a {\em curve} (resp. a {\em surface}) 
if $X$ is a variety of dimension one (resp. two).  
\item 
Let $f: X \to Y$ be a morphism of noetherian schemes. 
We say that $f$ is {\em projective} 
if it is projective in the sense of EGA \cite[D\'efinition 5.5.2]{EGAII}. 
This definition differs from the one in \cite[page 103]{Har77}. 
On the other hand, these coincide in many cases, e.g. $Y$ is affine (cf. \cite[Section 5.5.1]{FGAex}). 
\item For the definition of algberaic groups, we refer to \cite[Definition 1.1]{Mil17}. 
In particular, an algebraic group $G$ is of finite type over $k$, and however not necessarily reduced. 
\item 
A {\em torus} $T$ is an algebraic group which is isomorphic to $\mathbb G_m^r$ 
for some $r \geq 0$. 
We say that a variety $X$ is a {\em toric} if 
there exists a non-empty open subscheme $T$ of $X$ such that $T$ is a torus 
and the left action of $T$ on $T$ itself extends to an action on $X$: $T \times X \to X$. 
In this case, $T$ is called a {\em torus of} $X$. 
\item 
Throughout this note, $N$ denotes a free $\Z$-module of finite rank. 
Set $N_{\R} := N \otimes_{\Z} \R$. 
For terminologies on toric varieties, we refer to \cite{CLS11}. 
Although \cite{CLS11} works over $\mathbb C$, 
almost all the statements work even over arbitrary fields. 
For example, if $\Sigma$ is a fan of $N_{\R}$, then $X_{\Sigma}$ is a normal toric variety over $k$. 
\end{enumerate}

\subsection{Toric morphisms}

The purpose of this subsection is to recall the definition and some results on toric morphisms. 
All the materials treated in this subsection are well known to experts.

We start by recalling standard terminologies and some constructions on tori and lattices used in  toric geometry \cite[\S 1.1]{CLS11}. 
Let $T$ be a torus $T$, i.e., $T$ is an algebraic group over a field $k$ satisfying 
$T \simeq \G_m^n = \Spec\,k[t_1, t_1^{-1}, ..., t_n, t_n^{-1}]$ for some $n \in \Z_{\geq 0}$. 
Then we obtain free $\Z$-modules $N := \Hom(\G_m, T) \simeq \Z^{\oplus n}$ and $M :=\Hom (T, \G_m) \simeq \Z^{\oplus n}$ equipped with the natural bilinear map $N \times M \to \Hom(\G_m, \G_m) \simeq \Z, (\varphi, \psi) \mapsto \psi \circ \varphi$, 
where $\Hom(-, -)$ denotes the $\Z$-module consisting of the homomorphisms of algebraic groups. 
The opposite construction is  as follows: given a free $\Z$-module $N$ of finite rank, 
we set $M := \Hom_{\Z}(N, \Z)$ and $T_N := \Spec\,k[N]$, 
where $k[N]$ denotes the group ring over $k$ associated with $N$. 




\begin{dfn}\label{d-def-F}
{\em The category $\mathcal F$ of fans} is defined as follows. 
\begin{enumerate}
\item 
An object of $\mathcal F$ is a pair $(N, \Sigma)$ 
consisting of a finitely generated free $\Z$-module $N$ and 
a fan $\Sigma$ of $N_{\R} =N \otimes_{\Z} \R$. 
\item  
A morphism $\varphi : (N, \Sigma) \to (N', \Sigma')$ is a $\Z$-module homomorphism $\varphi: N \to N'$ such that 
given $\sigma \in \Sigma$, there exists $\sigma' \in \Sigma'$ such that $\varphi_{\R}(\sigma) \subset \sigma'$, 
where $\varphi_{\R} : N \otimes_{\Z} \R \to N' \otimes_{\Z} \R$ denotes the induced $\R$-linear map. 
\end{enumerate}
\end{dfn}

\begin{dfn}\label{d-def-T}
{\em The category $\mathcal T$ of toric varieties} is defined as follows. 
\begin{enumerate}
\item 
An object of $\mathcal T$ is a pair $(T, X)$, where 
$X$ is a variety and $T$ is a non-empty  open subscheme of $X$ such that 
$T$ is a torus and the left action $T \times T \to T$ of $T$ on $T$ extends to an action of $T$ on $X$: $T \times X \to X$. 
\item 
A morphism $f : (T, X) \to (T', X')$ is a $k$-morphism $f:X \to X'$ of $k$-schemes such that 
$f(T) \subset T'$ and the induced $k$-morphism $f|_T : T \to T'$ is a homomorphism of algebraic groups. 
\end{enumerate}

{\em The category $\mathcal T^{{\rm nor}}$ of normal toric varieties} 
is the full subcategory of $\mathcal T$ whose objects are the pairs $(T, X)$ such that $X$ is normal. 
We say that $(T, X)$ is a {\em toric variety} if $(T, X)$ is an object of $\mathcal T$. 
In this case, $T$ is called a {\em torus of} $X$. 
We say that $X$ is a {\em toric variety} if $(T, X)$ is an object of $\mathcal T$ for some torus $T \subset X$. 

We say that $f: (T_1, X_1) \to (T_2, X_2)$ is a {\em toric morphism} if 
both $(T_1, X_1)$ and $(T_2, X_2)$ are objects of $\mathcal T$ and 
$f$ is a morphism of $\mathcal T$. 
We say that a morphism $f: X_1 \to X_2$ of $k$-schemes is a {\em toric morphism} if 
there exist open subschemes $T_1 \subset X_1$ and $T_2 \subset X_2$ such that 
both $(T_1, X_1)$ and $(T_2, X_2)$ are objects of $\mathcal T$ and 
$f:(T_1, X_1) \to (T_2, X_2)$ is a morphism of $\mathcal T$. 
\end{dfn}

\begin{rem}
In \cite{CLS11}, Cox--Little--Schenck implicitly fix a torus $T$ of a normal toric variety $X_{\Sigma}$. 
\end{rem}

\begin{rem}\label{r-T-equiv}
Note that if $(T, X)$ is an object of $\mathcal T$, then 
an action $\alpha:T \times_k X \to X$ extended from the left action of $T$ on $T$ is unique. 
Indeed, two extended actions $\alpha, \beta : T \times_k X \to X$ coincide on the dense open subset $T \times_k T$, 
which implies that $\alpha = \beta$. 
Furthermore, if $f:(T, X) \to (T', X')$ is a morphism of $\mathcal T$, 
then the induced actions $\alpha:T \times_k X \to X$ and $\alpha':T' \times_k X' \to X'$ commute with given arrows: 
\[
\begin{CD}
T \times_k X @>\alpha >> X\\
@VVf_T \times_k fV @VVf V\\
T' \times_k X' @>\alpha' >> X'.\\
\end{CD}
\]
\end{rem}

\begin{thm}\label{t-F-T-equiv}
The functor 
\[
\eta: \mathcal F \to \mathcal T^{{\rm nor}}, \qquad (N, \Sigma) \mapsto (T_N, X_{\Sigma})
\]
is an equivalence of categories. 
\end{thm}

\begin{proof}[Sketch]
As this result is well known to experts, 
we only give a sketch of the proof. 
It follows from \cite[Corollary 3.1.8]{CLS11} that $\eta$ is essentially surjective. 
For each $i \in \{1, 2\}$, 
we take an object $(N_i, \Sigma_i) \in \mathcal F$ and set 
$(T_i, X_i) := \eta(N_i, \Sigma_i) = (T_{N_i}, X_{\Sigma_i})$. 
It is enough to show that 
\[
\widetilde{\eta}:{\rm Hom}_{\mathcal F}((N_1, \Sigma_1), (N_2, \Sigma_2)) 
\to 
{\rm Hom}_{\mathcal T}((T_{1}, X_{1}), (T_{2}, X_{2})) 
\]
is bijective. 
The injectivity of $\widetilde{\eta}$ is clear. 

Let us overview a proof of the surjectivity of $\widetilde{\eta}$. 
Fix a morphism $f : (T_{1}, X_{1}) \to (T_{2}, X_{2})$ of $\mathcal T$. 
For each $i \in \{1, 2\}$, 
let $\alpha_i : T_{i} \times_k X_{i} \to X_{i}$ be the action extended from the $T_{i}$-action on itself. 
By Remark \ref{r-T-equiv}, $f$ and $f_T := f|_{T_{1}} :T_{1} \to T_{2}$ commute with $\alpha_1$ and $\alpha_2$. 
Since we have $N_i \simeq \Hom(\mathbb G_m, T_i)$ for the $\Z$-module  $\Hom(\mathbb G_m, T_i)$ consisting of all the homomorphisms 
$\mathbb G_m \to T_i$ of algebraic groups, 
$f_T$ induces a $\Z$-module homomorphism $\varphi : N_1 \to N_2$. 
Fix $\sigma_1 \in \Sigma_1$. 
For its distinguished point $\gamma_{\sigma_1} \in X_1$ and its $T_1$-orbit $O(\sigma_1) := T_1 \cdot \gamma_{\sigma_1} \subset X_1$, 
we have that 
\[
f(O(\gamma_{\sigma_1})) = f(T_1 \cdot \gamma_{\sigma_1})= T_1 \cdot f(\gamma_{\sigma_1})
\subset T_2 \cdot f(\gamma_{\sigma_1}) = O(\gamma_{\sigma_2})
\]
for some $\sigma_2 \in \Sigma_2$ and its distinguished point $\gamma_{\sigma_2} \in X_2$. 
Using \cite[Proposition 3.2.2]{CLS11}, 
we can check that $\varphi_{\R}(\Relint(\sigma_1) \cap N_1)\subset \sigma_2$ 
for $\varphi_{\R} := \varphi \otimes_{\Z} \R$, 
where $\Relint(\sigma_1)$ denotes the relative interior of $\sigma_1$. 
By a purely convex geometric argument, 
we can show that 
$\varphi_{\R}(\Relint(\sigma_1) \cap N_1)\subset \sigma_2$ implies $\varphi_{\R}(\sigma_1) \subset \sigma_2$. 
Then $\varphi: (N_1, \Sigma_1) \to (N_2, \Sigma_2)$ is a morphism of $\mathcal F$ such that 
$\widetilde{\eta}(\varphi) = f$. 
\end{proof}

\subsection{Contraction images of toric varieties}

\begin{lem}\label{l-cont-image}
Let $H$ be an algebraic group and let $Y$ be a variety on which $H$ acts. 
Assume that there exists a non-empty open subset $V$ of $Y$ such that 
$V$ is $H$-stable and the induced $H$-action on $V$ is trivial. 
Then the $H$-action on $Y$ is trivial. 
\end{lem}

\begin{proof}
Recall that $Y^H$ (resp. $V^H$) is the largest closed subscheme of $Y$ (resp. $V$) 
such that $Y^H$ (resp. $V^H$) 
is $H$-stable and the induced $H$-action is trivial \cite[Chapter 7, b]{Mil17}. 
By \cite[Theorem 7.1]{Mil17}, 
we obtain $V^H \to Y^H$ which completes the following commutative diagram: 
\[
\begin{tikzcd}
V^H \arrow[rdd, bend right, "\simeq"'] \arrow[rrd, bend left] \arrow[rd, dashed]\\
& Y^H \cap V \arrow[r, hook] \arrow[d, hook] & Y^H \arrow[d, hook]\\
& V \arrow[r, hook] & Y. 
\end{tikzcd}
\]
Since the $V$-action on $V$ is trivial, we have $V = V^H$, and hence $V=Y^H \cap V =V^H$. 
Then the closed subscheme $Y^H$ contains a non-empty open subset $V$ of $Y$. 
Since $Y$ is an integral scheme, we get $Y = Y^H$, i.e., the $H$-action on $Y$ is trivial.     
\end{proof}

\begin{prop}\label{p-cont-image}
Let $f:X \to Y$ be a proper morphism of normal varieties with $f_*\MO_X =\MO_Y$. 
Assume that $X$ is a toric variety. 
Let $T_X$ be a torus of $X$ in the sense of Definition \ref{d-def-T}. 
Set $T_Y := f(T_X)$. 
Then $(T_Y, Y)$ is a toric variety such that $f:(T_X, X) \to (T_Y, Y)$ is a toric morphism. 
In particular, $f:X \to Y$ is a toric morphism. 
\end{prop}

The following  proof is due to Michel Brion. 

\begin{proof}
By $f_*\MO_X = \MO_Y$ and  \cite[Theorem 7.2.1]{Bri17}, 
there is a $T_X$-action on $Y$ such that $f$ is $T_X$-equivariant. 
Fix a $k$-rational point $x \in T_X(k) \subset  X(k)$ and set $y := f(x) \in Y(k)$. 
Since its orbit $O_y$, {which is equal to $f(O_x)$,} contains a non-empty open subset of $Y$ {(note that Chevalley's theorem implies that $f(O_x)$ is a constructible dense subset of $Y$)}, 
$O_y$ itself is a non-empty open subset of $Y$, as $O_y$ is a locally closed subset of $Y$ \cite[Proposition 1.65(b)]{Mil17}. 
Let $H$ be the stabiliser (isotropy) subgroup of $T_X$ at $y$, which is an algebraic subgroup of $T_X$. 
Since $O_y$ is a non-empty $H$-stable open subset of $Y$ whose $H$-action is trivial, 
the $H$-action on $Y$ is trivial (Lemma \ref{l-cont-image}). 
Therefore, the torus $T_Y := T_X/H$ acts on $Y$ and $T_Y = T_X/H \simeq O_y$ 
\cite[Corollary 7.13 and Proposition 7.17]{Mil17}. 
\qedhere


\end{proof}

In Example \ref{e-irrat} and Example \ref{e-char2}, 
we provide two examples satisfying the following property: 
$f:X \to Y$ is a finite surjective morphism of projective normal varieties such that $X$ is toric but $Y$ is not. 

\begin{ex}\label{e-irrat}
Let $k$ be an algebraically closed field. 
Then it is known that there is a counterexample to Noether's problem over $k$. 
In other words, there exists a finite group $G$ such that $k \subset k(G)$ is not a purely transcendental extension, where 
$k(G) := k(\{ x_g\}_{g \in G})^G$ for the natural $G$-action: $h \cdot x_g := x_{hg}$ \cite[Theorem 3.6]{Sal84}. 
In particular, this $G$-action extends to the corresponding affine space
\[
\mathbb A_k^{|G|} = \Spec\,k[ \{x_g\}_{g \in G}]. 
\]
This action further extends to the following compactification: 
\[
(\mathbb P^1)^{|G|} = \prod_{g \in G} \Proj\,k[x_g, y_g], 
\]
where the $G$-action is given by 
\[
h \cdot ( [a_g : b_g])_{g \in G} := ([a_{hg} : b_{hg}])_{g \in G} \qquad {\rm for} \qquad a_g, b_g \in k. 
\]
Therefore, $(\mathbb P^1)^{|G|}/G$ is not a rational variety, and hence not toric. 
\end{ex}

\begin{ex}\label{e-char2}
Let $k$ be an algebraically closed field of characteristic two. 
Fix a positive integer $m$.  
By \cite[Page 142]{Eke87}, 
there exists a finite purely inseparable morphism 
\[
\pi : Y \to \mathbb P^{2m+1}
\]
of degree two, where $Y$ is a smooth quadric hypersurface in $\mathbb P^{2m+2}$. 
The absolute Frobenius morphism $F:\mathbb P^{2m+1} \to \mathbb P^{2m+1}$ factors through $\pi$: 
\[
F : \mathbb P^{2m+1} \xrightarrow{f} Y \xrightarrow{\pi} \mathbb P^{2m+1}. 
\]
Then $f$ is a finite surjective purely inseparable morphism of degree $p^{2m}$. 
Recall that we have ${\rm Pic} (Y) \simeq {\rm Pic}(\mathbb P^{2m+2}) \simeq \Z$, 
where the first isomorphism holds by the Grothendieck-Lefschetz hyperplane section theorem \cite[Example X 2.2 and Th\'{e}or\`{e}me XI 3.18]{SGA2}. 
In particular, $\rho(Y)=1$.

In what follows, let us show that $Y$ is not toric. 
Otherwise, $Y$ is a smooth projective toric variety with $\rho(Y)=1$, which implies that $Y$ is a projective space. 
However, a quadric hypersurface $Y$ is not isomorphic to a projective space 
by $(-K_Y)^{d} = 2d^d \neq  (d+1)^d = (-K_{\mathbb P^d})^d$ for $d := \dim Y \geq 2$. 
\end{ex}

\subsection{Toric minimal model program}

In this subsection, we summarise notation and results on the toric minimal model program. 
In (\ref{d-toric-MMP}), we recall terminologies on toric minimal model program. 
In the proof of our main theorem, we need a description of toric flips in Proposition \ref{p-flip-blowup} 
(cf. the proof of Theorem \ref{t-toric-Qfac}). 
This result is known if $X$ is semi-projective \cite[Lemma 15.5.7]{CLS11}. 
Since we work with a slightly generalised setting, we give a proof of it for the sake of completeness. 
We start by giving the following characterisation of the quasi-projectivity for toric varieties with convex supports. 

\begin{lem}\label{l-toric-q-proj}
Let $\Sigma$ be a fan of $N_{\R}$ and let $X_{\Sigma}$ be the toric variety of $\Sigma$. 
Then the following are equivalent. 
\begin{enumerate}
\item $X_{\Sigma}$ is quasi-projective over $k$ and $|\Sigma|$ is a convex subset of $N_{\R}$. 
\item There exist a semi-projective normal toric variety $X'$ and a non-negative integer $r$ such that $X_{\Sigma} \simeq X' \times \mathbb G_m^r$. 
\item There exists a projective morphism $X_{\Sigma} \to Z$ to an affine variety $Z$. 
\item There exists a toric projective morphism $X_{\Sigma} \to Z$ to a toric normal affine variety $Z$. 
\end{enumerate}
\end{lem}

\begin{proof}
Let us show that (1) implies (2). 
Assume (1). 
Let $V_{\Sigma}$ be the smallest $\R$-vector subspace of $N_{\R}$ that contains $\Sigma$. 
Set $\overline{N} := N \cap V_{\Sigma}$ and set $\overline{\Sigma} :=\Sigma$, where 
we consider $\overline{\Sigma}$ as a fan in $\overline{N}_{\R}$. 
Note that $\overline{N}$ is spanned by a subset of a basis of $N$. 
It follows from \cite[Proposition 3.3.11]{CLS11} that 
\[
X_{\Sigma} \simeq X_{\overline{\Sigma}} \times \mathbb G_m^r. 
\]
Since $X_{\overline{\Sigma}}$ is quasi-projective over $k$ and $\overline{\Sigma}$ has full dimensional convex support in $\overline{N}_{\R}$, 
$X_{\overline{\Sigma}}$ is semi-projective \cite[Proposition 7.2.9]{CLS11}. 
Hence, (2) holds. 

It is clear that (2) implies (3). 
Let us show that (3) implies (4). 
Assume (3). 
Let $X_{\Sigma} \xrightarrow{g} W \to Z$ be the Stein factorisation of $X \to Z$, 
i.e. $W = \Spec\,\Gamma(X_{\Sigma}, \MO_{X_{\Sigma}})$. 
By $g_*\MO_{X_{\Sigma}} = \MO_W$, 
both $W$ and $g$ are toric (Proposition \ref{p-cont-image}). 
Thus, (4) holds. 

Let us show that (4) implies (1). 
Assume (4). 
It is clear that $X_{\Sigma}$ is quasi-projective over $k$. 
We have the $\Z$-module homomorphism $\varphi : N \to N'$  corresponding to  $X_{\Sigma} \to Z =: X_{\Sigma'}$. 
Since $X_{\Sigma} \to X_{\Sigma'}$ is proper, 
we have $\varphi_{\R}^{-1}(|\Sigma'|) = |\Sigma|$ 
for $\varphi_{\R} : N_{\R} \to N'_{\R}$ \cite[Theorem 3.4.11]{CLS11}.  
As $Z = X_{\Sigma'}$ is affine, $\Sigma'$ is convex. 
Take $\alpha, \beta \in |\Sigma| = \varphi_{\R}^{-1}(|\Sigma'|)$. 
Then $\varphi_{\R}(\alpha) \in |\Sigma'|$ and $\varphi_{\R}(\beta) \in |\Sigma'|$ imply 
\[
\varphi_{\R}(\alpha + \beta)  = \varphi_{\R}(\alpha) +\varphi_{\R}(\beta)  \in |\Sigma'|. 
\]
Therefore, it holds that $\alpha + \beta \in \varphi_{\R}^{-1}(|\Sigma'|) = |\Sigma|$. 
Thus, $|\Sigma|$ is convex, i.e. (1) holds. 
\end{proof}

\begin{rem}\label{r-torus-factor}
Let $S$ be an affine normal toric variety. 
Let $f : X_{\Sigma} \to X_{\Sigma'}$ be a toric projective birational morphism of 
toric normal varieties which are projective over $S$. 
In this case, the associated $\Z$-module homomorphism $\varphi: N \to N'$ of their lattices is an isomorphism. 
Hence, we may assume that  $\Sigma$ and $\Sigma'$ are fans of a common lattice $N$. 
Since $f$ is proper, we have that $|\Sigma| = |\Sigma'|$. 
By the proof of Lemma \ref{l-toric-q-proj}, there exists a toric projective birational morphism of semi-projective normal toric varieties 
\[
\overline f : X_{\overline{\Sigma}} \to X_{\overline{\Sigma}'}
\]
such that $f = \overline f \times \mathbb G_m^r$ for some $r \in \Z_{\geq 0}$. 
\end{rem}


\begin{nothing}[Terminologies on toric minimal model program]\label{d-toric-MMP}
Let $\alpha:X_{\Sigma} \to S$ be a projective toric morphism of normal toric varieties. 
Assume that $X_{\Sigma}$ is $\Q$-factorial and $S$ is affine. 
Let $D$ be an  $\R$-divisor on $X_{\Sigma}$.

Note that $\NE(X_{\Sigma}/S)$ is generated by the proper curves on $X_{\Sigma}$ and 
$\NE(X_{\Sigma}/S)$ is a polyhedral cone \cite[Theorem 4.1]{FS04}, 
i.e. $\NE(X_{\Sigma}/S) = \sum_{i=1}^r \R_{\geq 0} [C_i]$ for some proper curves $C_1, ..., C_r$ on $X_{\Sigma}$. 
Furthermore, we may and do assume that each $ \R_{\geq 0} [C_i]$ is an extremal ray of $\NE(X_{\Sigma}/S)$. 
Let $R =\R_{\geq 0}[C_i]$ be an extremal ray of $\NE(X_{\Sigma}/S)$. 
Then there exists {\em the contraction of} $R$ \cite[Theorem 4.5]{FS04}, i.e. the projective toric morphism 
\[
\varphi_R : X_{\Sigma} \to X_{\Sigma'}
\]
to a normal toric variety $X_{\Sigma'}$ projective over $S$ 
such that $(\varphi_R)_*\MO_{X_{\Sigma}} = \MO_{X_{\Sigma'}}$ and 
that if $C$ is a proper curve on $X_{\Sigma}$, then $[C] \in R$ if and only if $\varphi_R(C)$ is a point. 
We say that $R =\R_{\geq 0}[C_i]$ is {\em $D$-negative} if $D \cdot C_i <0$. 
\begin{itemize}
\item We say that $\varphi_R$ is a {\em $D$-divisorial contraction} if $R$ is $D$-negative, 
$\varphi_R$ is birational, and $\dim \Ex(\varphi_R) =\dim X_{\Sigma}-1$. 
In this case, $X_{\Sigma'}$ is $\Q$-factorial \cite[Theorem 4.5]{FS04}. 
\item We say that $\varphi_R$ is a {\em $D$-flipping contraction} if $R$ is $D$-negative,  $\varphi_R$ is birational, and $\dim \Ex(\varphi_R) <\dim X_{\Sigma}-1$. 
\item We say that $\varphi_R$ is a {\em $D$-Mori fibre space} if $R$ is $D$-negative and  $\varphi_R$ is not birational, i.e. $\dim X_{\Sigma} > \dim X_{\Sigma'}$. 
In this case, also $X_{\Sigma}$ is called a $D$-Mori fibre space. 
\end{itemize}
If $\varphi_R : X_{\Sigma} \to X_{\Sigma'}$ is a $D$-flipping contraction, then there exists a small projective toric morphism 
$\varphi^+ : X_{\Sigma^+} \to X_{\Sigma'}$ from a $\Q$-factorial toric variety $X_{\Sigma^+}$ such that the proper transform $D^+$ of $D$ on  $X_{\Sigma^+}$ is $\varphi^+$-ample \cite[Theorem 4.8]{FS04}. 
In this case, each of $X_{\Sigma} \dashrightarrow  X_{\Sigma^+}$, $\varphi^+$, and $X_{\Sigma^+}$ is called a 
{\em $D$-flip} of $\varphi_R$. 

By the same argument as in \cite[Theorem 4.9]{FS04}, there exists a sequence
\[
X_{\Sigma} = X_{\Sigma_0} \overset{f_0}{\dashrightarrow} X_{\Sigma_1}  \overset{f_1}{\dashrightarrow} \cdots 
 \overset{f_{\ell-1}}{\dashrightarrow} X_{\Sigma_{\ell}}, 
\]
called a {\em $D$-MMP over $S$}, such that the following hold for $D_{i+1} := (f_i)_* D_i$. 
\begin{enumerate}
\item For any $i \in \{0, ..., \ell\}$, $\alpha_i : X_i \to S$ is a toric projective morphism 
from a $\Q$-factorial toric variety $X_i$. 
\item For any $i \in \{0, ..., \ell-1\}$, 
$f_i: X_{\Sigma_i} \dashrightarrow X_{\Sigma_{i+1}}$ is either a $D_i$-divisorial contraction over $S$ 
or a $D_i$-flip over $S$. 
\item 
Either $D_{\ell}$ is nef over $S$ or there exists a $D_{\ell}$-Mori fibre space $X_{\Sigma_{\ell}} \to X_{\Lambda}$ 
over $S$. 
\end{enumerate}
\end{nothing}

\begin{rem}\label{r-MMP-torus-factor}
Let $\alpha:X_{\Sigma} \to S$ be a projective toric morphism of normal toric varieties. 
Assume that $X_{\Sigma}$ is $\Q$-factorial and $S$ is affine. 
Let $D$ be an $\R$-divisor on $X_{\Sigma}$. 
\begin{enumerate}
\item 
We have $X_{\Sigma} \simeq X_{\overline{\Sigma}} \times \mathbb G_m^r$ for some $r \in \Z_{\geq 0}$ and 
a $\Q$-factorial semi-projective toric variety $X_{\overline{\Sigma}}$. 
By ${\rm Pic}(X_{\overline{\Sigma}} \times \mathbb G_m^r) \simeq {\rm Pic}(X_{\overline{\Sigma}})$, 
there exists an $\R$-divisor $\overline D$ on $X_{\overline{\Sigma}}$ 
such that $\pi^* \overline{D} \sim_{\R} D$, where $\pi: X_{\Sigma} \to X_{\overline{\Sigma}}$ denotes the natural projection. 
\item 
Let $\varphi : X_{\Sigma} \to X_{\Sigma'}$ be the contraction of an extremal ray $R$ of $\NE(X_{\Sigma}/S)$. 
Assume that $\varphi$ is birational. 
Then we have $\varphi  = \overline{\varphi} \times \mathbb G_m^r$ for a suitable projective birational morphism 
$\overline{\varphi} : X_{\overline \Sigma} \to X_{\overline{\Sigma}'}$ of semi-projective toric normal varieties 
(Remark \ref{r-torus-factor}). 
In this case, we have $\rho(X_{\overline \Sigma}/X_{\overline{\Sigma}'})=1$ and $-\overline{D}$ is $\overline \varphi$-ample, and 
hence $\overline{\varphi}$ is the contraction of some extremal ray $\overline{R}$ of $\NE(X_{\Sigma})$. 
Furthermore, $\varphi$ is a $D$-divisorial contraction (resp. $D$-flipping contraction) if and only if 
$\overline \varphi$ is a $\overline D$-divisorial contraction (resp. $\overline D$-flipping contraction). 
\item 
Let $\varphi : X_{\Sigma} \to X_{\Sigma'}$ be a $D$-divisorial contraction. 
Then $E := \Ex(\varphi)$ is a prime divisor on $X_{\Sigma}$. 
Furthermore, $E$ is torus-invariant. 
Indeed, this is known if $X_{\Sigma}$ is semi-projective \cite[Proposition 15.4.5(b)]{CLS11}, and hence the general case is reduced to this case by (2). 
\item 
Let $\varphi : X_{\Sigma} \to X_{\Sigma'}$ be a $D$-flipping contraction and let $\varphi^+ : X_{\Sigma^+} \to X_{\Sigma'}$ 
be its flip. 
By Remark \ref{r-torus-factor}, 
we have a flip $\overline{\varphi}^+ : X_{\overline{\Sigma}^+} \to X_{\overline{\Sigma}'}$ of $\overline{\varphi} : X_{\overline \Sigma} \to X_{\overline{\Sigma}'}$ such that $\varphi^+ = \overline{\varphi}^+ \times \mathbb G_m^r$. 
\end{enumerate}
\end{rem}

\begin{prop}\label{p-flip-blowup}
Let $\alpha:X_{\Sigma} \to S$ be a projective toric morphism of normal toric varieties. 
Assume that $X_{\Sigma}$ is $\Q$-factorial and $S$ is affine. 
Let $D$ be an $\R$-divisor on $X_{\Sigma}$. 
Let $X_{\Sigma} \to X_{\Pi}$ be a $D$-flipping contraction and let $X_{\Sigma'} \to X_{\Pi}$ be its flip. 
Then there exists a commutative diagram consisting of a projective toric birational $S$-morphisms of toric normal varieties which are projective over $S$ 
\[
\begin{tikzcd}
	& X_{\Theta} \\
	{X_{\Sigma}} && {X_{\Sigma'}} \\
	& X_{\Pi}
	\arrow["\psi"', from=1-2, to=2-1]
	\arrow["\varphi"', from=2-1, to=3-2]
	\arrow["{\psi'}", from=1-2, to=2-3]
	\arrow["{\varphi'}", from=2-3, to=3-2]
	\arrow["{f}", dashed, from=2-1, to=2-3]
\end{tikzcd}
\]
such that the following hold. 
\begin{enumerate}
\item $X_{\Theta}$ is $\Q$-factorial. 
\item $E := \Ex(\psi) = \Ex(\psi')$ is a torus-invariant prime divisor on $X_{\Theta}$. 
\item There exists a nonzero effective $1$-cycle $\Gamma$ on $X_{\Sigma}$ such that 
\begin{equation}\label{e1-flip-blowup}
\psi^* F = \psi'^* F'  - (F \cdot \Gamma) E
\end{equation}
holds for any $\R$-divsior $F$ on $X_{\Sigma}$ and $F' := f_*F$. 
\end{enumerate}
\end{prop}

\begin{proof}
By Remark \ref{r-MMP-torus-factor}(2), 
we also have $\overline{\varphi} : X_{\overline{\Sigma}} \to X_{\overline{\Pi}}$
and 
$\overline{\varphi}' : X_{\overline{\Sigma}'} \to X_{\overline{\Pi}}$ such that 
$\varphi = \overline{\varphi} \times \mathbb G_m^r$ and $\varphi' = \overline{\varphi}' \times \mathbb G_m^r$. 
Since $\overline{\Sigma}$ is semi-projective, 
it follows from \cite[Lemma 15.5.7]{CLS11} that 
we have the following right square diagram. 
Then, by taking the direct product $(-) \times \mathbb G_m^r$, 
we obtain the following left square diagram (cf. Remark \ref{r-torus-factor}, e.g. $X_{\Theta} = X_{\overline{\Theta}} \times \mathbb G_m^r$). 
\[
\begin{tikzcd}
	& X_{\Theta} \\
	{X_{\Sigma}} && {X_{\Sigma'}} \\
	& X_{\Pi}
	\arrow["\psi"', from=1-2, to=2-1]
	\arrow["\varphi"', from=2-1, to=3-2]
	\arrow["{\psi'}", from=1-2, to=2-3]
	\arrow["{\varphi'}", from=2-3, to=3-2]
	\arrow["{f}", dashed, from=2-1, to=2-3]
\end{tikzcd}
\hspace{10mm}
\begin{tikzcd}
	& X_{\overline{\Theta}} \\
	{X_{\overline{\Sigma}}} && {X_{\overline{\Sigma}'}} \\
	& X_{\overline{\Pi}}
	\arrow["\overline{\psi}"', from=1-2, to=2-1]
	\arrow["\overline{\varphi}"', from=2-1, to=3-2]
	\arrow["{\overline{\psi}'}", from=1-2, to=2-3]
	\arrow["{\overline{\varphi}'}", from=2-3, to=3-2]
	\arrow["{\overline{f}}", dashed, from=2-1, to=2-3]
\end{tikzcd}
\]

Let us show (1). 
Since $X_{\overline{\Sigma}}$ and $X_{\overline{\Sigma}'}$ are $\Q$-factorial, 
$X_{\overline \Theta}$ is $\Q$-factorial by \cite[Lemma 15.5.7(c)]{CLS11}. 
Then $\overline{\Theta}$ is simplicial, and hence also $\Theta$ is simplicial. 
Therefore, $X_{\Theta}$ is $\Q$-factorial, i.e. (1) holds. 

Let us show (2). 
By \cite[Lemma 15.5.7(c)]{CLS11}, there exists a unique $\psi$-exceptional prime divisor $\widetilde E$ on $X_{\Theta}$. 
Pick an effective ample Cartier divisor $H$ on $X_{\Theta}$ whose support does not contain any irreducible component of $\Ex(\psi)$. 
Since 
$X_{\Sigma}$ is $\Q$-factorial, $\psi_*H$ is $\Q$-Cartier. 
Therefore, $\psi^*\psi_*H - H$ is a $\psi$-exceptional effective divisor, and hence we have $ \psi^*\psi_*H - H = a\widetilde{E}$ for some $a \in \Q_{\geq 0}$. 
Since $a\widetilde{E}$ is $\psi$-anti-ample, we have $a>0$ and $\widetilde E \cdot C <0$ for any curve $C$ on $X_{\Theta}$ such that 
$C \subset \Ex(\psi)$. This implies $\Ex(\psi) \subset \widetilde E$. The opposite inclusion is clear, and hence $\widetilde E =\Ex(\psi)$. 
The same argument implies $\widetilde{E}=\Ex(\psi')$.

Let us show (3). The equation (\ref{e1-flip-blowup}) holds if $\Sigma$ is semi-projective \cite[Lemma 15.5.7(c)]{CLS11}. 
Hence, we have a nonzero effective $1$-cycle $\overline \Gamma$ on $X_{\overline \Sigma}$ such that 
\[
\overline{\psi}^* \overline{F} = \overline{\psi}^* \overline{F}'  - (\overline{F} \cdot \overline{\Gamma}) \overline{E}
\]
holds for any $\Q$-divsior $\overline F$ on $X_{\overline \Sigma}$ and $\overline F' := \overline f_*\overline F$. 
For the closed point $\{1\} \in \mathbb G_m^r$, 
let $\Gamma$ be the $1$-cycle on $X \times \{ 1\} \subset X_{\Sigma} \times \mathbb G_m^r$ 
that is the pullback of $\overline \Gamma$ via the composite isomorphism $X \times \{ 1\} \hookrightarrow  X_{\Sigma} \times \mathbb G_m^r \to X_{\Sigma}$. 

Let us show the equation (\ref{e1-flip-blowup}) for any $\R$-divisor $F$ on $X_{\Sigma}$. 
By linearity, the problem is reduced to the case when $F$ is a Cartier divisor. 
In order to prove the equation (\ref{e1-flip-blowup}),  
we may replace $F$ by a Cartier divisor linearly equivalent to $F$. 
Therefore, it suffices to prove the equation (\ref{e1-flip-blowup}) 
for the case when $F = \pi^*\overline{F}$ for some Cartier divisor $\overline{F}$ on $X_{\overline{\Sigma}}$ (cf. Remark \ref{r-MMP-torus-factor}(1)). 
We then obtain 
 \[
\psi^* F =\psi^*\pi^*\overline F = \psi'^* \pi'^*\overline F'- (\overline F \cdot \overline \Gamma) E
= \psi^* F'  - (F \cdot \Gamma) E, 
\]
where $\pi:X_{\Sigma} \to X_{\overline{\Sigma}}$ and $\pi':X_{\Sigma'} \to X_{\overline{\Sigma}'}$ denote 
the natural projections. 
Thus (3) holds. 
\end{proof}

\subsection{Globally $F$-regular varieties}

In this subsection, we first recall the definition of globally $F$-regular pairs 
(Definition \ref{d-GFR}) and a relative version (Definition \ref{d-rel-GFR}). 
Strictly speaking, we do not need the relative version in our proof. 
However, this will make our presentation clearer. 
We also give a proof of the folklore fact that $(X_{\Sigma}, (1-\epsilon)D)$ is globally $F$-regular, 
where $X_{\Sigma}$ is a normal toric variety and $D$ is the sum of the torus invariant prime divisors 
(Proposition \ref{p-toric-GFR}). 

\begin{dfn}[cf. Definition 3.1 of \cite{SS10}]\label{d-GFR}
Assume that $k$ is of characteristic $p>0$. 
Let $X$ be a normal variety and let $\Delta$ be an effective $\R$-divisor on $X$. 
\begin{enumerate}
\item 
We say that $(X, \Delta)$ is {\em sharply globally $F$-split} 
if 
there exists $e \in \Z_{>0}$ such that 
\[
\MO_X \xrightarrow{F^e} F_*^e\MO_X \hookrightarrow F_*^e(\MO_X(\ulcorner (p^e-1)\Delta \urcorner))
\]
splits as an $\MO_X$-module homomorphism. 
\item 
We say that $(X, \Delta)$ is {\em globally $F$-regular} 
if {\cred for any} 
effective $\Z$-divisor $E$, there exists $e \in \Z_{>0}$ such that 
\[
\MO_X \xrightarrow{F^e} F_*^e\MO_X \hookrightarrow F_*^e(\MO_X((p^e-1)\Delta + E))
\]
splits as an $\MO_X$-module homomorphism. 
\end{enumerate}
\end{dfn}

\begin{dfn}[cf. Definition 2.6 of \cite{HX15}]\label{d-rel-GFR}
Assume that $k$ is of characteristic $p>0$. 
Let $f:X \to Y$ be a morphism from a normal variety $X$ to a scheme $Y$ of finite type over $k$. 
Let $\Delta$ be an effective $\R$-divisor on $X$. 
\begin{enumerate}
\item 
We say that $(X, \Delta)$ is {\em sharply globally $F$-split over $Y$} 
if 
there exists an open cover $Y = \bigcup_{i \in I} Y_i$ such that 
$(f^{-1}(Y_i), \Delta|_{f^{-1}(Y_i)})$ is sharply globally $F$-split for any $i \in I$. 
\item 
We say that $(X, \Delta)$ is {\em globally $F$-regular over $Y$} 
if 
there exists an open cover $Y = \bigcup_{i \in I} Y_i$ such that 
$(f^{-1}(Y_i), \Delta|_{f^{-1}(Y_i)})$ is globally $F$-regular for any $i \in I$. 
\end{enumerate}
\end{dfn}

\begin{prop}\label{p-toric-GFR}
Assume that $k$ is of characteristic $p>0$. 
Let $\Sigma$ be a fan of $N_{\R}$ and let $X_{\Sigma}$ be the toric variety of $\Sigma$. 
Set $D := \sum_{\rho \in \Sigma(1)} D_{\rho}$ 
(recall that $D_{\rho}$ denotes the prime divisor that is the closure of the orbit $O(\rho)$ of $\rho \in \Sigma(1)$). 
Then the following hold. 
\begin{enumerate}
\item $(X, D)$ is globally sharply $F$-split. 
\item $(X, (1-\epsilon)D)$ is globally $F$-regular for any $\epsilon \in \R$ with $0 < \epsilon \leq 1$. 
\end{enumerate}
\end{prop}

\begin{proof}
Set $X := X_{\Sigma}$ {\cred and let  $\overline k$ be the algebraic closure of  $k$}. 
{\cred Since $X$ is toric, 
also the base change $X_{\overline k}$ is toric  
(e.g., if $X$ is affine, i.e., $X = \Spec\,k[S_{\sigma}]$ for $S_{\sigma} := M \cap \sigma^{\vee}$, then we have $k[S_{\sigma}] \otimes_k \overline k = \overline{k}[S_{\sigma}]$). In particular, $X_{\overline k}$ is a normal variety.}   
{\cred Moreover,} 
if the base change $(X_{\overline k}, \Delta_{\overline k})$ is globally sharply $F$-split (resp. globally $F$-regular), 
then so is $(X, \Delta)$, 
{\cred because $\MO_X \to \theta_*\MO_{X_{\overline k}}$ splits for the induced morphism $\theta : X_{\overline k} \to X$}. 
Hence, the problem is reduced to the case when $k$ is an algebraically closed field. 

Let us show (1). 
For the relative Frobenius $\overline{F}: X \to X'$ and the absolute Frobenius morphisms $F:X \to X$ and $F_k : \Spec\,k \to \Spec\,k$, 
we have the following commutative diagram in which the square diagram is cartesian: 
\[
\begin{tikzcd}
	X & {X'} & X \\
	& \Spec\,k & \Spec\,k
	\arrow["{\overline F}", from=1-1, to=1-2]
	\arrow["\alpha", from=1-2, to=1-3]
	\arrow["{F_k}", from=2-2, to=2-3]
	\arrow["\pi", from=1-3, to=2-3]
	\arrow["{\pi'}", from=1-2, to=2-2]
	\arrow["{F =F_X}", bend left, from=1-1, to=1-3]
	\arrow["\pi"', from=1-1, to=2-2]. 
\end{tikzcd}
\]
Note that $\alpha$ and $F_k$ are isomorphisms of schemes because $k$ is a perfect field. 
For $D' := \alpha^*D$, it follows from \cite[2.6]{Fuj07} that 
\[
\MO_X(-D) \to \overline F_*\MO_{X'}(-D'). 
\]
splits {\cred (note that $\Omega_X^a(\log D)= \MO_X(-D)$ if $a =0$)}. 
Since $\alpha$ is an isomorphism, also 
\[
F : \MO_X(-D) \to F_*\MO_X(-D)
\]
splits. Applying $(-) \otimes_{\MO_X} \MO_X(D)$, we see that 
\[
F : \MO_X \to F_*\MO_X( (p-1)D)
\]
splits. Thus, (1) holds.

Let us show (2). 
Fix a rational number $\epsilon$ with $0 < \epsilon <1$. 
Set $\Delta := (1-\epsilon)D$. 
By \cite[Theorem 3.9]{SS10}, it suffices to show the following assertions (i) and (ii). 
\begin{enumerate}
\renewcommand{\labelenumi}{(\roman{enumi})}
\item There exists $e \in \Z_{>0}$ such that 
\[
\MO_X \xrightarrow{F^e} F_*^e\MO_X \hookrightarrow F_*^e\MO_X( \ulcorner (p^e-1) \Delta \urcorner + D)
\]
splits as an $\MO_X$-module homomorphism. 
\item 
$(X \setminus \Supp\,D, \Delta|_{X\setminus \Supp\,D})$ is globally $F$-regular. 
\end{enumerate}

Let us show (i). 
By $\Delta = (1-\epsilon)D$ and $0 < \epsilon <1$, 
we have 
\[
(p^e-1) \Delta + D \leq (p^e-1) D
\]
for $e \gg 0$. Since 
\[
\MO_X \to F_*^e\MO_X( (p^e-1)D)
\]
splits by (1), we see that (i) holds. 
The assertion (ii) follows from the fact that $\Delta|_{X\setminus \Supp\,D} =0$ and 
$X \setminus \Supp\,D$ is isomorphic to the torus $\mathbb G_m^{\dim X}$ 
(note that $(Z, 0)$ is globally $F$-regular for any smooth affine variety $Z$). 
\end{proof}

\section{Vanishing theorems for globally F-regular varieties}

In this section, we establish some vanishing theorems for globally $F$-regular varieties. 
In particular, we establish an analogue of our main theorem (Theorem \ref{t-intro-main1}) 
for globally $F$-regular surfaces (Theorem \ref{t-GFR-surface}). 

\begin{lem}\label{l-rel-minus-big}
Let $\alpha :X \to S$ be a proper morphism of normal varieties 
such that $\alpha_*\MO_X=\MO_S$ and $\dim X> \dim S$. 
Let $D$ be a $\Q$-Cartier $\Z$-divisor on $X$ such that $-D$ is $\alpha$-big. 
Then $\alpha_*\MO_X(D)=0$. 
\end{lem}

\begin{proof}
We may assume that $S$ is affine. 
It is enough to show that $H^0(X, \MO_X(D))=0$. 
Suppose that $H^0(X, \MO_X(D)) \neq 0$. 
Then there exists an effective $\Q$-Cartier $\Z$-divisor $E$ on $X$ such that 
$D \sim E$. 
For the generic fibre $X_{\xi}$ of $\alpha:X \to S$, 
we have that $E|_{X_{\xi}}$ is effective and $-E_{X_{\xi}}$ is big. 
This contradicts $\dim X_{\xi} =\dim X - \dim S>0$. 
\end{proof}

Although the following proposition is well known to experts, 
we give a proof for the sake of completeness. 

\begin{prop}\label{p-GFR-easy}
Assume that $k$ is of characteristic $p>0$. 
Let $\alpha:X \to S$ be a projective morphism 
from a normal variety $X$ to a scheme $S$ of finite type over $k$. 
Let $\Delta$ be an effective $\R$-divisor on $X$ such that $(X, \Delta)$ is globally $F$-regular over $S$. 
\begin{enumerate}
\item 
If $D$ is a $\Q$-Cartier $\Z$-divisor on $X$ such that $D-(K_X+\Delta)$ is $\alpha$-nef and $\alpha$-big, 
then $R^if_*\MO_X(D) = 0$ for $i>0$. 
\item 
If $D$ is an $\alpha$-nef $\Q$-Cartier $\Z$-divisor on $X$, then 
$R^if_*\MO_X(D)=0$ for any $i>0$. 
\end{enumerate}
\end{prop}

\begin{proof}
Taking the Stein factorisation of $S$ and an affine cover of $S$, 
we may assume that $\alpha_*\MO_X = \MO_S$, $S$ is affine, and $(X, \Delta)$ is globally $F$-regular. 
Furthermore, the problem is reduced to the case when $k$ is $F$-finite by taking a suitable subfield of $k$ which is finitely generated over $\mathbb F_p$. 

Let us show (1). 
Since $D-(K_X+\Delta)$ is $\alpha$-nef and $\alpha$-big, 
there exist $\epsilon_1 \in \Q_{>0}$  and an effective $\Q$-Cartier $\Q$-divisor $E$ such that 
\[
D-(K_X+\Delta) -\epsilon E
\]
is $\alpha$-ample for any $\epsilon \in \Q$ with $0< \epsilon <\epsilon_1$. 
There also exists $\epsilon_2 \in \Q_{>0}$ such that $(X, \Delta +\epsilon E)$ is globally $F$-regular 
if $0 \leq \epsilon < \epsilon_2$ \cite[Corollary 6.1]{SS10}. 
After replacing $(X, \Delta)$ by $(X, \Delta + \epsilon E)$ for $0 < \epsilon < \min \{\epsilon_1, \epsilon_2\}$, 
we may assume that $D-(K_X+\Delta)$ is $\alpha$-ample.

We have a split surjection {\cred (cf. \cite[Section 4.2]{SS10})}: 
\[
F_*^e\MO_X( -\ulcorner (p^e-1)(K_X+\Delta)\urcorner) \to \MO_X. 
\]
Taking the tensor product with $\MO_X(D)$ and applying the double dual, we obtain another split surjection: 
\[
F_*^e\MO_X( p^eD-\ulcorner (p^e-1)(K_X+\Delta)\urcorner) \to \MO_X(D). 
\]
For $i>0$ and $e \gg 0$, we have 
\begin{eqnarray*}
&& H^i(X, F_*^e\MO_X( p^eD-\ulcorner (p^e-1)(K_X+\Delta)\urcorner) )\\
&\simeq& H^i(X, \MO_X( p^eD-\ulcorner (p^e-1)(K_X+\Delta)\urcorner) )\\
&=& H^i(X, \MO_X( D + \llcorner (p^e-1)(D-(K_X+\Delta))\lrcorner) )\\
&=& 0, 
\end{eqnarray*}
where the last equality follows from the Serre vanishing theorem. 
Therefore, it holds that 
\[
H^i(X, \MO_X(D))=0, 
\]
which completes the proof of (1). 

Let us show (2). 
Since $(X, \Delta)$ is globally $F$-regular, 
after possibly 
{\cred replacing} 
$S$ by a suitable open cover, there exists 
an effective $\R$-divisor  $\Delta'$ on $X$ such that 
$-(K_X+\Delta)$ is $\alpha$-ample \cite[Theorem 4.3(i)]{SS10}. 
In particular, $D-(K_X+\Delta')$ is $\alpha$-ample, hence 
the assertion (2) follows from (1). 
\end{proof}

\begin{prop}\label{p-GFR-main}
Assume that $k$ is of characteristic $p>0$. 
Let $\alpha :X \xrightarrow{f} Y \xrightarrow{\beta} S$ be projective morphisms of normal varieties 
with $f_*\MO_X = \MO_Y$. 
Let $\Delta$ be an effective $\R$-divisor 
such that $K_X+\Delta$ is $\R$-Cartier and $(X, \Delta)$ is globally $F$-regular over $Y$. 
Let $D$ be a $\Q$-Cartier $\Z$-divisor such that 
\begin{enumerate}
\item $-D$ is $f$-big, 
\item $D -(K_X+\Delta)$ is $f$-nef and $f$-big, and 
\item $\dim X>\dim Y$. 
\end{enumerate}
Then $R^i\alpha_*\MO_X(D)=0$ for any $i \geq 0$. 
\end{prop}

\begin{proof}
We have the following isomorphism in the derived category: 
\[
R\alpha_*\MO_X(D) \simeq 
R\beta_* Rf_*\MO_X(D). 
\]
Hence, it suffices to show that $Rf_*\MO_X(D)=0$, i.e. $R^if_*\MO_X(D) =0$ for any $i \geq 0$. 
We have $f_*\MO_X(D)=0$ by (1), (3), and Lemma \ref{l-rel-minus-big}. 
By (2), it follows from Proposition \ref{p-GFR-easy} that  $R^if_*\MO_X(D) =0$ holds for any $i >0$. 
\end{proof}

\begin{cor}\label{c-GFR-main}
Let $\alpha :X \to S$ be a projective morphism from a normal $\Q$-factorial variety $X$ to a scheme $S$ of finite type over $k$. 
Let $D$ be a $\Z$-divisor on $X$ such that $D -K_X$ is $\alpha$-big. 
Assume that 
\begin{enumerate}
\item $X$ is globally $F$-regular over $S$, and 
\item either $D$ is $\alpha$-nef or 
there exists a $D$-Mori fibre space $f:X \to Y$ over $S$, i.e. 
$f$ is a projective $S$-morphism to a normal variety $Y$ projective over $S$ 
to a normal variety $Y$ which is projective over $S$ such that $f_*\MO_X = \MO_Y$, 
$\dim X > \dim Y$, $\rho(X/Y)=1$, and $-D$ is $f$-ample.  
\end{enumerate}
Then $R^i\alpha_*\MO_X(D)=0$ for any $i>0$. 
\end{cor}

\begin{proof}
If $D$ is $\alpha$-nef, then the assertion follows from Proposition \ref{p-GFR-easy}(2). 
Hence, we may assume that there exists a $D$-Mori fibre space $f:X \to Y$ over $S$. 
Then we have the induced morphisms $\alpha:X \xrightarrow{f} Y \xrightarrow{\beta} S$. 
Since $D -K_X$ is $\alpha$-big, $D-K_X$ is $f$-big. 
Furthermore, it follows from $\dim X > \dim Y$ and $\rho(X/Y)=1$ that any effective $\Q$-divisor on $X$ is $f$-nef. 
Therefore, $D-K_X$ is $f$-nef and $f$-big. 
Since $-D$ is $f$-ample, it holds by Proposition \ref{p-GFR-main} that 
$R^i\alpha_*\MO_X(D)=0$ for any $i>0$. 
\end{proof}


\begin{thm}\label{t-GFR-surface}
Assume that $k$ is of characteristic $p>0$. 
Let $\alpha :X \to S$ be a proper morphism 
from a normal surface $X$ to a scheme $S$ of finite type over $k$ 
such that  $X$ is globally $F$-regular over $S$. 
Let $D$ be a 
$\Z$-divisor on $X$ such that 
$D -(K_X+B)$ is $\alpha$-nef and $\alpha$-big 
for some $\R$-divisor $B$ on $X$ whose coefficients are contained in $[0, 1)$. 
Then $R^i\alpha_*\MO_X(D)=0$ for any $i > 0$. 
\end{thm}

\begin{proof}
Note that $X$ is $\Q$-factorial \cite[Corollary 4.11]{Tan18}. 
By the same argument as in \cite[Lemma 2.2]{Fuj12}, $X$ is quasi-projective over $k$, 
and hence $f$ is projective. 
Taking a suitable affine open cover of $S$, 
we may assume that $S$ is affine and $X$ is globally $F$-regular. 
We run a $D$-MMP over $S$ \cite[Theorem 1.1]{Tan18}: 
\[
X = X_0 \xrightarrow{f_0} X_1 \xrightarrow{f_1} \cdots \xrightarrow{f_{N-1}} X_N. 
\]
Set $D_n$ and $B_n$ to be the push-forwards of $D$ and $B$ to $X_n$, respectively. 
Note that either $D_N$ is $\alpha_N$-nef or there exists a $D_N$-Mori fibre space $X_N \to Y$ over $S$. 
{\cred Moreover, each $X_i$ is globally $F$-regular over $S$ \cite[Proposition 1.2(2)]{HWY02}.} 
By Corollary \ref{c-GFR-main}, we get 
\[
H^i(X_N, \MO_{X_N}(D_N))=0
\]
for any $i>0$. 
Fix $n \in \{0, 1, ..., N-1\}$. 
It is enough to show that 
\[
H^i(X_n, \MO_{X_n}(D_n)) \simeq H^{i}(X_{n+1}, \MO_{X_{n+1}}(D_{n+1})).
\]
For the $f_n$-exceptional prime divisor $E := \Ex(f_n)$, 
we have $D_n = f_n^*D_{n+1} + a E$ for some $a \in \Q_{>0}$. 
Therefore, we obtain 
\[
(f_n)_*\MO_{X_n}(D_n) = \MO_{X_{n+1}}(D_{n+1}). 
\]
We have the following spectral sequence:  
\[
E_2^{i, j} =H^i(X_{n+1}, R^j(f_n)_*\MO_{X_{n}}(D_n)) \Rightarrow H^{i+j}(X_n, \MO_{X_n}(D_n)). 
\]
By \cite[Theorem 3.3]{Tan18}, 
we obtain $R^j(f_n)_*\MO_{X_{n}}(D_n) =0$ for $j>0$. 
Hence, 
\[
H^i(X_{n+1}, \MO_{X_{n+1}}(D_{n+1}))
=
H^i(X_{n+1}, (f_n)_*\MO_{X_{n}}(D_n)) 
\]
\[
=E_2^{i, 0} \simeq E^i = H^{i}(X_n, \MO_{X_n}(D_n)), 
\]
which completes the proof. 
\end{proof}

\section{Vanishing theorems for toric varieties}

In this section, we prove the main theorem of this paper (Theorem \ref{t-main-full}). 
To this end, we first treat a special but essential case (Theorem \ref{t-toric-Qfac}). 
In order to conclude Theorem \ref{t-main-full} from Theorem \ref{t-toric-Qfac}, 
we also establish an auxiliary result (Lemma \ref{l-toric-nef}).

\begin{thm}\label{t-toric-Qfac}
Let $\alpha: X \to S$ be a projective morphism 
from a $\Q$-factorial normal toric variety $X$ to a separated scheme $S$ of finite type over $k$. 
Let $D$ be a $\Z$-divisor on $X$ such that 
$D \sim_{\R} K_X+B$ for some klt pair $(X, B)$, where $B$ is an $\alpha$-big effective $\R$-divisor.  
Then $R^i\alpha_*\MO_X(D)=0$ for any $i>0$. 
\end{thm}

\begin{proof}
We may assume that $k$ is of characteristic $p>0$, as otherwise the assertion is well known. 
We divide the proof into four steps. 
\setcounter{step}{0}

\begin{step}\label{s1-toric-Qfac}
In order to prove Theorem \ref{t-toric-Qfac}, we may assume that the following conditions hold. 
\begin{enumerate}
\renewcommand{\labelenumi}{(\roman{enumi})}
\item $\alpha_*\MO_X = \MO_S$. 
\item $S$ is an affine normal toric variety. 
\item $f$ is a toric morphism. 
\end{enumerate}
\end{step}

\begin{proof}(of Step \ref{s1-toric-Qfac}) 
Taking the Stein factorisation of $f$, we may assume (i). 
Then the problem is reduced to the case when $f$ is a toric morphism to a normal toric variety $S$ 
(Proposition \ref{p-cont-image}).
In particular, (iii) holds. 
Replacing $(X, S)$ by $(f^{-1}(S_i), S_i)$ for a toric open affine cover $S = \bigcup_{i \in I} S_i$, 
we may further assume that (ii) holds. 
This completes the proof of Step \ref{s1-toric-Qfac}. 
\end{proof}

In what follows, we assume that (i)--(iii) hold. 
We run a $D$-MMP (i.e. a $(K_X+B)$-MMP) over $S$ (cf. (\ref{d-toric-MMP})): 
\[
X = X_0 \overset{f_0}{\dashrightarrow} X_1 \overset{f_1}{\dashrightarrow} \cdots  \overset{f_{N-1}}{\dashrightarrow} X_N
\]
\[
B_{n+1} := (f_n)_*B_n, \qquad D_{n+1} := (f_n)_*D_n, 
\]
where $X_N$ denotes the end result, so that 
either $D_N$ is nef over $S$ or there exists a $D_N$-Mori fibre space $X_N \to W$ over $S$. 
Note that the assumptions are stable under this MMP, i.e. 
for any $n \in \{0, 1, ..., N\}$, it holds that 
\begin{enumerate}
\renewcommand{\labelenumi}{(\roman{enumi})}
\item[(iv)] the induced morphism $\alpha_n : X_n \to S$ is a projective morphism, 
\item[(v)] $X_n$ is a $\Q$-factorial toric variety, 
\item[(vi)] $(X_n, B_n)$ is klt, 
\item[(vii)] $B_n$ is an $\alpha_n$-big effective $\R$-divisor, and 
\item[(viii)] $D_n$ is a $\Z$-divisor such that $D_n \sim_{\R} K_{X_n} +B_n$. 
\end{enumerate}


\begin{step}\label{s2-toric-Qfac}
The equation $H^i(X_N, \MO_{X_N}(D_N)) =0$ holds for any $i>0$.  
\end{step}

\begin{proof}(of Step \ref{s2-toric-Qfac}) 
Since ${\cred X_N}$ is a toric variety, ${\cred X_N}$ is globally $F$-regular (Proposition \ref{p-toric-GFR}). 
In particular, ${\cred X_N}$ is globally $F$-regular over $S$. 
It follows from Corollary \ref{c-GFR-main} that $H^i({\cred X_N}, \MO_{{\cred X_N}}({\cred D_N}))=0$ for any $i>0$. 
This completes the proof of Step \ref{s2-toric-Qfac}. 
\end{proof}

\begin{step}\label{s3-toric-Qfac}
Fix $n \in \{0, 1, ..., N-1\}$. 
Assume that $f_n : X_n \dashrightarrow X_{n+1}$ is a $D_n$-divisorial contraction. 
Then 
\[
H^i(X_n, \MO_{X_n}(D_n)) \simeq H^i(X_{n+1}, \MO_{X_{n+1}}(D_{n+1}))
\]
for any $i \geq 0$. 
\end{step}

\begin{proof}(of Step \ref{s3-toric-Qfac}) 
Since $f_n : X_n \to X_{n+1}$ is a  $D_n$-divisorial contraction, 
$E:=\Ex(f_n)$ is a torus invariant primer divisor (Remark \ref{r-MMP-torus-factor}(3)). 
Then we have 
\[
D_n =f_n^*D_{n+1} + a E
\]
for some $a \in \Q_{>0}$. 
In particular, we obtain $(f_n)_*\MO_{X_n}(D_n) = \MO_{X_{n+1}}(D_{n+1})$. 
Let us prove that $R^i(f_n)_*\MO_{X_n}(D_n)=0$ for any $i>0$. 
For the coefficient $b$ of $E$ in $B_n$, 
we can write $B_n = bE + A$, where $A$ is $f_n$-nef and $f_n$-big. 
We have that 
\[
D \sim_{\R} K_{X_n} + B_n = K_{X_n} + b E+A_n
\]
Since $(X_n, bE)$ is globally $F$-regular (Proposition \ref{p-toric-GFR}), we have that $R^i(f_n)_*\MO_{X_n}(D)=0$ for $i>0$ (Proposition \ref{p-GFR-easy}(1)). 
Therefore it holds that 
\[
H^i(X_n, \MO_{X_n}(D_n)) \simeq H^i(X_{n+1}, \MO_{X_{n+1}}(D_{n+1})). 
\]
This completes the proof of Step \ref{s3-toric-Qfac}. 
\end{proof}

\begin{step}\label{s4-toric-Qfac}
Fix $n \in \{0, 1, ..., N-1\}$. 
Assume that $f_n : X_n \dashrightarrow X_{n+1}$ is a $D_n$-flip. 
Then 
\[
H^i(X_n, \MO_{X_n}(D_n)) \simeq H^i(X_{n+1}, \MO_{X_{n+1}}(D_{n+1}))
\]
for any $i \geq 0$. 
\end{step}

\begin{proof}(of Step \ref{s4-toric-Qfac}) 
We have the $D_n$-flipping contraction $\varphi: X_n \to Z$ and its flip $\varphi' : X_{n+1} \to Z$. 
In this case, by Proposition \ref{p-flip-blowup}, 
we have the following square diagram of projective birational toric morphisms of quasi-projective normal toric varieties  
\[
\begin{tikzcd}
	& Y \\
	{X_n} && {X_{n+1}} \\
	& Z
	\arrow["\psi"', from=1-2, to=2-1]
	\arrow["\varphi"', from=2-1, to=3-2]
	\arrow["{\psi'}", from=1-2, to=2-3]
	\arrow["{\varphi'}", from=2-3, to=3-2]
	\arrow["{f_n}", dashed, from=2-1, to=2-3]
\end{tikzcd}
\]
such that 
\begin{enumerate}
\item $Y$ is  $\Q$-factorial, 
\item $E:=\Ex(\psi) = \Ex(\psi')$ is a torus invariant prime divisor on $Y$, and 
\item there exists a nonzero effective $1$-cycle $\Gamma$ on $X_n$ satisfying the equation 
\[
\psi^* F = \psi'^*( (f_n)_*F) -(F \cdot \Gamma) E 
\]
for any $\R$-divisor $F$ on $X_n$. 
\end{enumerate}
We have 
\begin{equation}\label{e1-s4-toric-Qfac}
\psi^*D_{n} \sim_{\R} \psi^*(K_{X_{n}} +B_{n})
\end{equation}
and 
\begin{equation}\label{e2-s4-toric-Qfac}
K_Y + \psi_*^{-1}B_{n} = \psi^*(K_{X_{n}}+B_{n}) + aE 
\end{equation}
for some rational number $a>-1$. 
By (\ref{e1-s4-toric-Qfac}) and (\ref{e2-s4-toric-Qfac}), we obtain   
\begin{equation}\label{e3-s4-toric-Qfac}
\psi^*D_{n} \sim_{\R} K_Y + \psi_*^{-1}B_{n} -a E. 
\end{equation}
On the other hand, we have 
\begin{equation}\label{e4-s4-toric-Qfac}
\ulcorner \psi^*D_{n} \urcorner = \psi^*D_{n} +bE
\end{equation} 
for some rational number $b \in [0, 1)$. 
Since $\varphi : X_n \to Z$ is the contraction of a $D_n$-negative extremal ray, 
it holds by (3) that 
\begin{equation}\label{e5-s4-toric-Qfac}
\psi^* D_n = \psi'^*D_{n+1} + c E 
\end{equation}
for some rational number $c > 0$. 
By (\ref{e3-s4-toric-Qfac}) and (\ref{e4-s4-toric-Qfac}), we get 
\begin{equation}\label{e6-s4-toric-Qfac}
\ulcorner \psi^*D_{n} \urcorner \sim_{\R} K_Y + \psi_*^{-1}B_{n} +(-a+b) E. 
\end{equation}
We separately treat the following two cases: $-a +b <1$ and $-a+b \geq 1$.

Assume that $-a +b <1$. Then there uniquely exists $m \in \Z_{\geq 0}$ such that $-a+b +m \in [0, 1)$. 
Set 
\[
D_Y := \ulcorner \psi^*D_{n} \urcorner +mE \sim_{\R} K_Y + \psi_*^{-1}B_{n} +(-a+b+m) E. 
\]
Since $E$ is torus invariant, $(Y, (-a+b+m)E)$ is globally $F$-regular (Proposition \ref{p-toric-GFR}). 
Note that $\psi_*^{-1}B_{n}$ is $\psi$-nef and $\psi'$-nef. 
It follows from Proposition \ref{p-GFR-easy}(1) that 
\[
R^i\psi_*\MO_Y(D_Y)= R^i\psi'_*\MO_Y(D_Y)=0
\]
for any $i>0$. Hence, we have 
\[
H^i(X_n, \psi_*\MO_Y(D_Y)) \simeq H^i(X_{n+1}, \psi'_*\MO_Y(D_Y)) 
\]
for any $i \geq 0$. 
By $c>0$ and 
\[
D_Y = \ulcorner \psi^*D_{n} \urcorner +mE = \ulcorner \psi'^*D_{n+1} + cE \urcorner +mE,  
\]
we have $\psi_*\MO_Y(D_Y ) = \MO_{X_n} (D_n)$ and $\psi'_*\MO_Y(D_Y) = \MO_{X_{n+1}}(D_{n+1})$, so that we get 
\[
H^i(X_n, \MO_{X_n}(D_n)) \simeq H^i(X_{n+1}, \MO_{X_{n+1}}(D_{n+1})) 
\]
for any $i \geq 0$. 
This completes the proof for the case when $-a +b <1$.

We may assume that $-a +b \geq 1$. 
In particular, we obtain $0 < -a< 1$ and $0<b <1$. 
Hence, $0 \leq -a +b -1 <1$. Set $D_Y := \llcorner \psi^*D_n \lrcorner$. 
We obtain 
\[
D_Y = \llcorner \psi^*D_n \lrcorner = \ulcorner \psi^*D_n \urcorner - E = \psi^*D_n -(1-b)E 
\]
\[
\sim_{\R} K_Y + \psi_*^{-1}B_n  +(-a+b-1)E. 
\]
where the second equality holds by (\ref{e4-s4-toric-Qfac}) and $0<b<1$, 
the third one by (\ref{e4-s4-toric-Qfac}), 
and the last $\R$-linear equivalence by (\ref{e3-s4-toric-Qfac}). 
Since $E$ is torus invariant, $(Y, (-a+b-1)E)$ is globally $F$-regular (Proposition \ref{p-toric-GFR}). 
As $\psi_*^{-1}B_n$ is $\psi$-nef and $\psi'$-nef, 
it follows from Proposition \ref{p-GFR-easy} that 
\[
R^j\psi_*\MO_Y(D_Y)= R^j\psi'_*\MO_Y(D_Y)=0
\]
for any $j>0$, which in turn implies 
\[
H^i(X_n, \psi_*\MO_Y(D_Y)) \simeq H^i(X_{n+1}, \psi'_*\MO_Y(D_Y)) 
\]
for any $i \geq 0$. 
By $D_Y = \llcorner \psi^*D_n \lrcorner = \llcorner \psi'^*D_{n+1} + cE \lrcorner$ and $c >0$, we have 
\[
\psi_*\MO_{Y}(D_Y) = \MO_{X_n}(D_n) \qquad {\rm and }\qquad \psi'_*\MO_Y(D_Y) =\MO_{X_{n+1}}(D_{n+1}). 
\]
Hence, we obtain 
\[
H^i(X_n, \MO_{X_n}(D_n)) \simeq H^i(X_{n+1}, \MO_{X_{n+1}}(D_{n+1})), 
\]
This completes the proof of Step \ref{s4-toric-Qfac}. 
\end{proof}
Step \ref{s2-toric-Qfac}, Step \ref{s3-toric-Qfac}, and Step \ref{s4-toric-Qfac} complete the proof of 
Theorem \ref{t-toric-Qfac}. 
\end{proof}

\begin{lem}\label{l-toric-nef}
Let $\alpha : X \to S$ be a proper morphism from a normal toric variety to a scheme of finite type over $k$. 
Let $D$ be an $\alpha$-nef $\R$-Cartier $\R$-divisor. 
Then $D$ is $\alpha$-semi-ample. 
\end{lem}

\begin{proof}
We may assume that $S$ is affine. 
By toric Chow's lemma, the problem is reduced to the case when $\alpha$ is projective. 
Since ${\rm NE}(X/S)$ is a rational polyhedtral cone \cite[Theorem 4.1]{FS04}, 
its dual ${\rm Nef}(X/S)$ is a rational polyhedral cone, 
where ${\rm Nef}(X/S)$ consists of the classes of the $\alpha$-nef $\R$-divisors on $X$. 
Then we can write $D = \sum_{i=1}^{r} a_i D_i$ for some $a_1, ..., a_r \in \R_{>0}$ and 
nef $\Q$-divisors $D_1, ..., D_r$. 
Since each $D_i$ is $\alpha$-semi-ample \cite[Proposition 4.6]{FS04}, also $D$ is $\alpha$-semi-ample 
(note that \cite[Proposition 4.6]{FS04} imposes the assumption that $D_i$ is torus invariant, 
however we can drop this assumption because any Cartier divisor is linearly equivalent to a torus invariant Cartier divisor \cite[Theorem 4.2.1]{CLS11}). 
\end{proof}

\begin{thm}\label{t-main-full}
Let $\alpha: X \to S$ be a projective morphism 
from a normal toric variety $X$ to a separated scheme $S$ of finite type over $k$. 
Let $D$ be a $\Q$-Cartier $\Z$-divisor on $X$. 
Assume that one of the following conditions holds. 
\begin{enumerate}
\item $D \sim_{\R} K_X+B_1+B_2$ for some effective $\R$-divisors $B_1$ and $B_2$ on $X$ 
such that $(X, B_1 +B_2)$ is klt and $B_2$ is an $\alpha$-big $\R$-Cartier $\R$-divisor. 
\item $D - (K_X+B)$ is $\alpha$-nef and $\alpha$-big for some effective $\R$-divisor $B$ on $X$ 
such that $(X, B)$ is klt. 
\end{enumerate}
Then $R^i\alpha_*\MO_X(D)=0$ for any $i>0$. 
\end{thm}

\begin{proof}
We may assume that $S$ is affine, $\alpha_*\MO_X = \MO_S$, and $\alpha$ is a toric morphism to a normal toric variety $S$ (Proposition \ref{p-cont-image}).

Assume (1). 
Let $\mu : X' \to X$ be a projective toric small birational morphism 
from a $\Q$-factorial toric variety $X'$ \cite[Proposition 11.1.7]{CLS11}:
\[
\alpha' : X' \xrightarrow{\mu} X \xrightarrow{\alpha} S. 
\]
We define a $\Z$-divisor $D'$ and $\R$-divisors $B'_1$ and  $B'_2$ on $X'$ by 
\[
\mu_* B'_1 = B_1, \quad 
\mu_* B'_2 = B_2, \quad 
\mu_* D' = D. 
\]
Note that $D'$ is a $\Z$-divisor. 
Set $B' := B'_1 + B'_2$, which is an effective $\R$-divisor on $X'$. 
We then have 
\begin{enumerate}
\item[(a)] $B'_2 = \mu^*B_2$, 
\item[(b)] $K_{X'} + B' = \mu^*(K_X+B_1 + B_2)$, and 
\item[(c)] $D' =\mu^*D \sim_{\R} K_{X'} + B'$. 
\end{enumerate}
Since $(X, B_1 + B_2)$ is klt, also $(X', B')$ is klt by (b). 
As $B_2$ is $\alpha$-big, $B'_2$ is $\alpha'$-big by (a). 
Hence, $B' =B'_1 + B'_2$ is $\alpha'$-big. 
By Theorem \ref{t-toric-Qfac}, we obtain 
\[
R^i\alpha'_*\MO_{X'}(D') = 0 
\]
for any $i>0$. 
Since  $B'$ is $\mu$-big, it follows again from Theorem \ref{t-toric-Qfac} that 
\[
R^j\mu_*\MO_{X'}(D') = 0 
\]
for any $j>0$. By the Leray spectral sequence 
\[
E_2^{i, j} := H^i(X, R^j\mu_*\MO_{X'}(D')) \Rightarrow H^{i+j}(X', \MO_{X'}(D')) =: E^{i+j}, 
\]
we obtain 
\[
H^i(X, \MO_{X}(D)) = H^i(X, \mu_*\MO_{X'}(D')) = E_2^{i, 0} \simeq E^i = H^{i}(X', \MO_{X'}(D')) = 0 
\]
for any $i>0$, where the first equality is guaranteed by $D' = \mu^*D$. 
This completes the proof for the case when (1) holds. 

Assume (2). It suffices to show that (1) holds. 
Set $A := D-(K_X+B)$. 
Note that $A$ is $\alpha$-semi-ample (Lemma \ref{l-toric-nef}). 
In particular, we may assume that $A$ is semi-ample, i.e. $|mA|$ is base point free for some $m \in \Z_{>0}$. 
It follows from \cite[Theorem 1.4]{Tanb} that after taking a suitable base change of the base field $k$, 
there exists an effective $\R$-Cartier $\R$-divisor $B_2$ such that $A \sim_{\R} B_2$ and $(X, B+B_2)$ is klt. 
We see that the condition (1) holds for $B_1 := B$. 
\end{proof}






\end{document}